\documentclass[12pt,a4paper]{article}
\usepackage{amsfonts,amssymb,graphicx}
\newtheorem{thm}{Theorem}
\newcounter{abbildung}
\newcommand{\cB}{{\mathcal B}}
\newcommand{\cM}{{\mathcal M}}
\newcommand{\cP}{{\mathcal P}}
\newcommand{\cS}{{\mathcal S}}
\newcommand{\cV}{{\mathcal V}}
\newcommand{\cW}{{\mathcal W}}
\newcommand{\qed}{\ $\square$}
\newcommand{\pf}{{\sl{Proof:}\ }}
\newcommand{\abstand}{\vspace{0.5em}}
\begin{document}
\sloppy

\date{}
\title{Another Simple Proof for the\\ Existence of the Small Witt Design}

\author{Hans Havlicek%
\thanks{Partially supported by the City of Vienna
({\em Hochschuljubil\"aumsstiftung der Stadt Wien, Projekt
H-39/98}.)}
\and {Hanfried Lenz}}

\maketitle

\begin{abstract} \noindent
We give a short proof for the existence of the small Witt
design which is based on the projective plane of order three
with one point deleted.
\end{abstract}

\section{Introduction}

The Swiss geometer J.~Steiner posed the following question
(``Combinatorische Aufgabe'') in 1853:

``Welche Zahl, $N$, von Elementen hat die Eigenschaft, dass sich die
Elemente so zu dreien ordnen lassen, dass je zwei in einer, aber nur
in einer Verbindung vorkommen?''

If we write $v$, $k$, and $t$ instead of $N$, $3$, and $2$,
respectively, then we arrive at the following contemporary
definition: A {\em Steiner system} $S(t,k,v)$ is a finite set $\cV$
of elements (called {\em points}) with a distinguished family of
subsets (called {\em blocks}) such that the following holds true:
\begin{enumerate}
  \item There are exactly $v$ points in $\cV$.
  \item Each block has exactly $k$ elements.
  \item Any $t$ distinct points belong to a unique block.
\end{enumerate}
In order to avoid trivialities it is usually assumed that $2\leq t<
k< v$.

So Steiner asked for $S(2,3,v)$ systems. As a matter of fact,
T.P.~Kirkman proved already in 1847 that an $S(2,3,v)$ exists if, and
only if, $v\equiv 1,\,3 \pmod 6$.

In this short communication we present another proof for the
existence of a Steiner system $S(5,6,12)$ which is also called {\em
small Witt design} $W_{12}$. See \cite{beth-jung-lenz85} or
\cite{beth-jung-lenz99}, in particular Chapter IV. There the reader
will also find the definition of a {\em $t$-design} (which is more
general than that of a Steiner system) and references on other
results mentioned in this section.

In an $S(5,6,12)$ there are twelve points, each block has exactly six
elements, and any five distinct points are contained in a unique
block. There is a unique $S(5,6,12)$ up to isomorphism. The same
uniqueness property holds true for an $S(5,8,24)$ which carries the
name {\em large Witt design $W_{24}$}. The Steiner systems $W_{12}$
and $W_{24}$ are due to E.~Witt (1938) and R.D.~Carmichel (1937). For
many decades $W_{12}$ and $W_{24}$ were the only known Steiner
systems with parameter $t=5$. Even today only finitely many Steiner
systems $S(t,k,v)$ with $t>3$ and none with $t>5$ seem to be known
\cite[67]{colb-math-96}, \cite{math97}.

Another remarkable property of the two Witt designs concerns their
automorphism group. Recall that a group $G$ of permutations acts
(sharply) $t$-transitively, if for two ordered $t$-tuples of elements
there is a (unique) permutation in $G$ taking the first to the second
$t$-tuple. The automorphism groups of the Witt designs $W_{12}$ and
$W_{24}$ act $5$-transitively on their sets of points; for the small
Witt design the action is even sharply $5$-transitive. These
automorphism groups are the {\em Mathieu groups} $M_{12}$ and
$M_{24}$, respectively. They were discovered by E.~Mathieu in 1861
and 1873, and they are early examples of {\em sporadic finite simple
groups}. The only finite $t$-transitive permutation groups with $t>3$
other than symmetric and alternating groups that seem to be known are
the two Mathieu groups mentioned above and two of their subgroups
(the Mathieu groups $M_{11}$ and $M_{23}$). So the Witt designs are
indeed remarkable combinatorial structures.

The starting point of our construction of $W_{12}$ is the projective
plane of order three with point set $\cP$. It is a Steiner system
$S(2,4,13)$, but its blocks are called {\em lines}.

The first step is to discuss the $6$-sets of points in $\cP$. They
fall into four classes which can be described in various ways, but
the crucial observation is that two $6$-sets are in the same class if
and only if they have the same number of trisecants (i.e.\ lines
meeting the set in exactly three points).

Next we choose one point of $\cP$, say $U$. The twelve
points of $\cW:=\cP\setminus\{U\}$ will be the points of the
Witt design $W_{12}$. We introduce three kinds of
$6$-subsets of $\cW$ and call them blocks. Each block
together with the distinguished point $U$ has a complement
in $\cP$ with exactly six elements. So properties of
$6$-sets in $\cP$ carry over to properties of blocks.

Finally, we show that $\cW$, together with the set of all
blocks, is a Steiner system $S(5,6,12)$. Again, the results
on $6$-sets of points turn out useful when showing that any
$5$-set $\cM\subset\cW$ is contained in a block, since
$\cM\cup\{U\}$ is a $6$-set of points in the projective
plane.

The proof presented in this paper is closely related to a
projective representation, in the five-dimensional
projective space of order three, of the small Witt design
due to H.S.M.\ Coxeter \cite{coxe58}; see \cite{havl99} and
the references given there. Furthermore, we refer to
\cite{havl0x} for an alternative description of the present
construction of $W_{12}$ using completely different methods.


\section{Construction}

Let $\cP$ be the set of points of the projective plane of order three
or, in other words, the Steiner system $S(2,4,13)$
\cite[19]{beth-jung-lenz85}. There are exactly $4$ lines (blocks)
through each point of $\cP$. The unique line joining distinct points
$A$ and $B$ will be written as $AB$.

First we introduce four types of sets $\cS\subset\cP$, each
consisting of exactly six points.

\begin{enumerate}
\item
$\cS$ is the union of a line and two further points (fig.\
\ref{abb1}).
\item
$\cS$ is the symmetric difference of two different lines
(fig.\ \ref{abb2}).
\item
$\cS$ consists of a triangle and an inscribed triangle,
i.e.\ each point of the second triangle lies on exactly one
line of the first triangle (fig.\ \ref{abb3}).
\item
$\cS$ is the set of vertices of a quadrilateral, i.e.\ the
set of points where two distinct lines of the quadrilateral
meet (fig.\ \ref{abb4}).
\end{enumerate}

{\unitlength1cm
      \begin{center}
      \begin{minipage}[t]{4.0cm}
         \begin{picture}(4.0,4.0)
         \put(0,0)
         {\includegraphics[width=4.0cm,height=4.0cm]{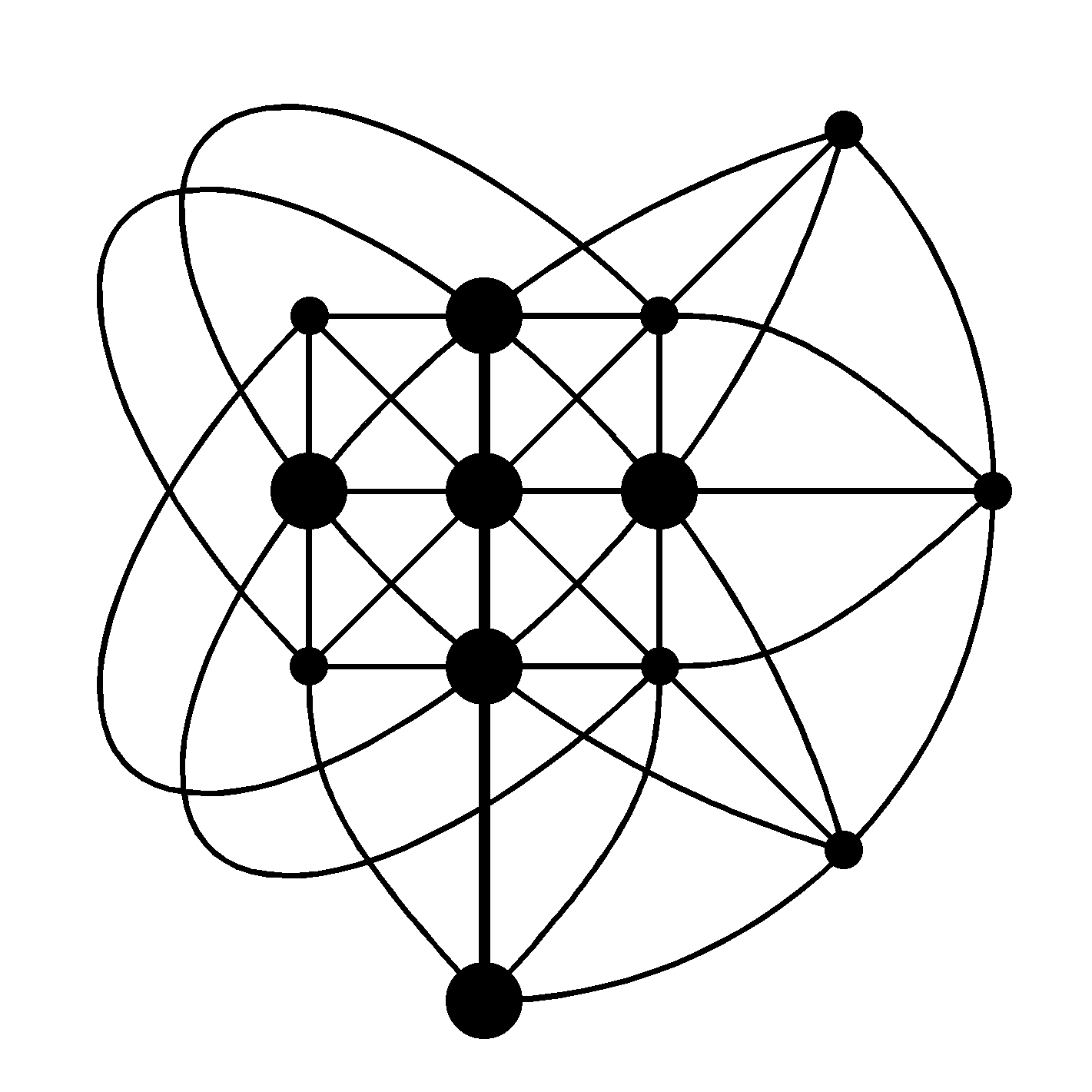}}
         \end{picture}
         {\refstepcounter{abbildung}\label{abb1}
          \centerline{Fig.\ \ref{abb1}.}}
         \end{minipage}
      \begin{minipage}[t]{4.0cm}
         \begin{picture}(4.0,4.0)
         \put(0,0)
         {\includegraphics[width=4.0cm,height=4.0cm]{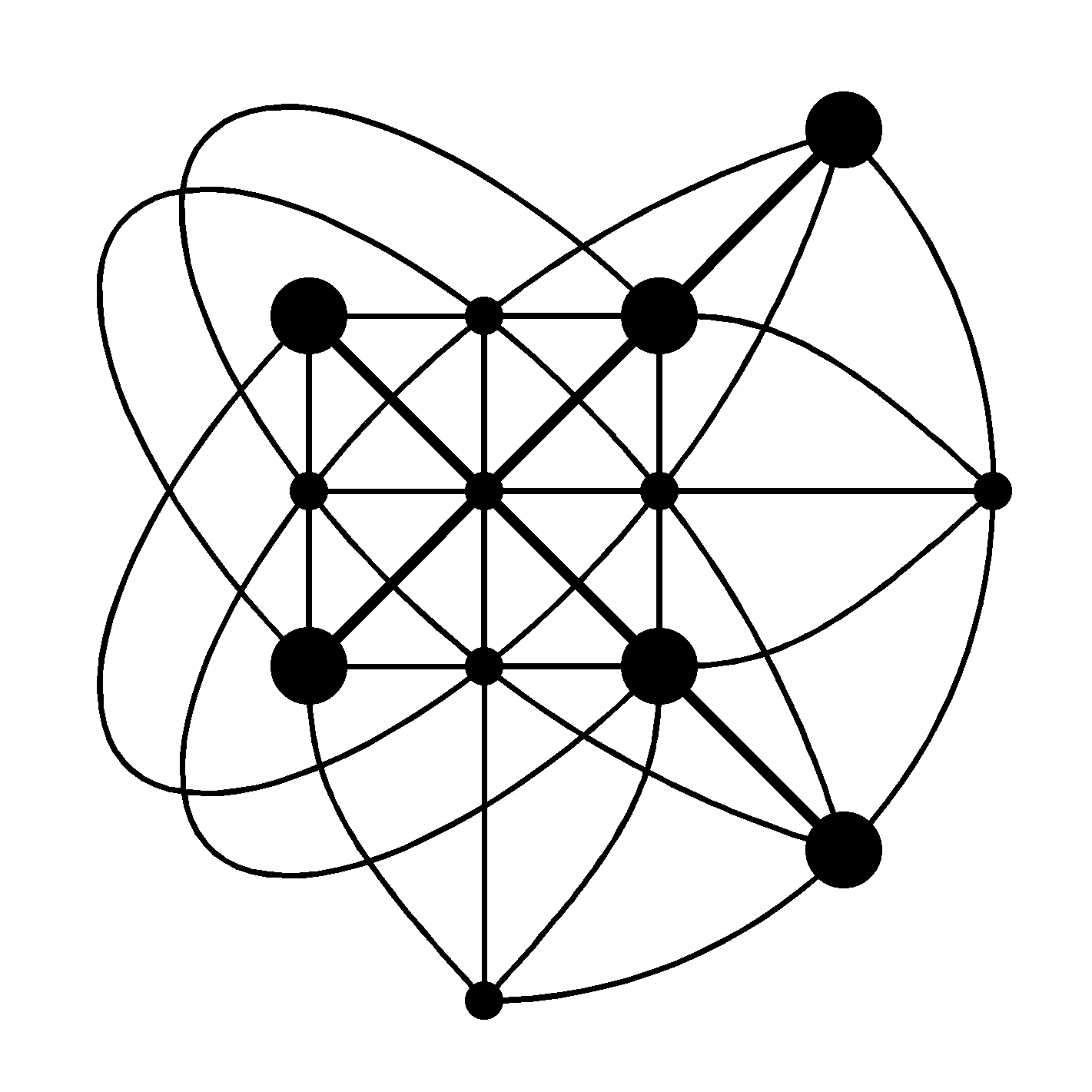}}
        \end{picture}
        {\refstepcounter{abbildung}\label{abb2}
          \centerline{Fig.\ \ref{abb2}.}}
      \end{minipage}
      \begin{minipage}[t]{4.0cm}
         \begin{picture}(4.0,4.0)
         \put(0,0)
         {\includegraphics[width=4.0cm,height=4.0cm]{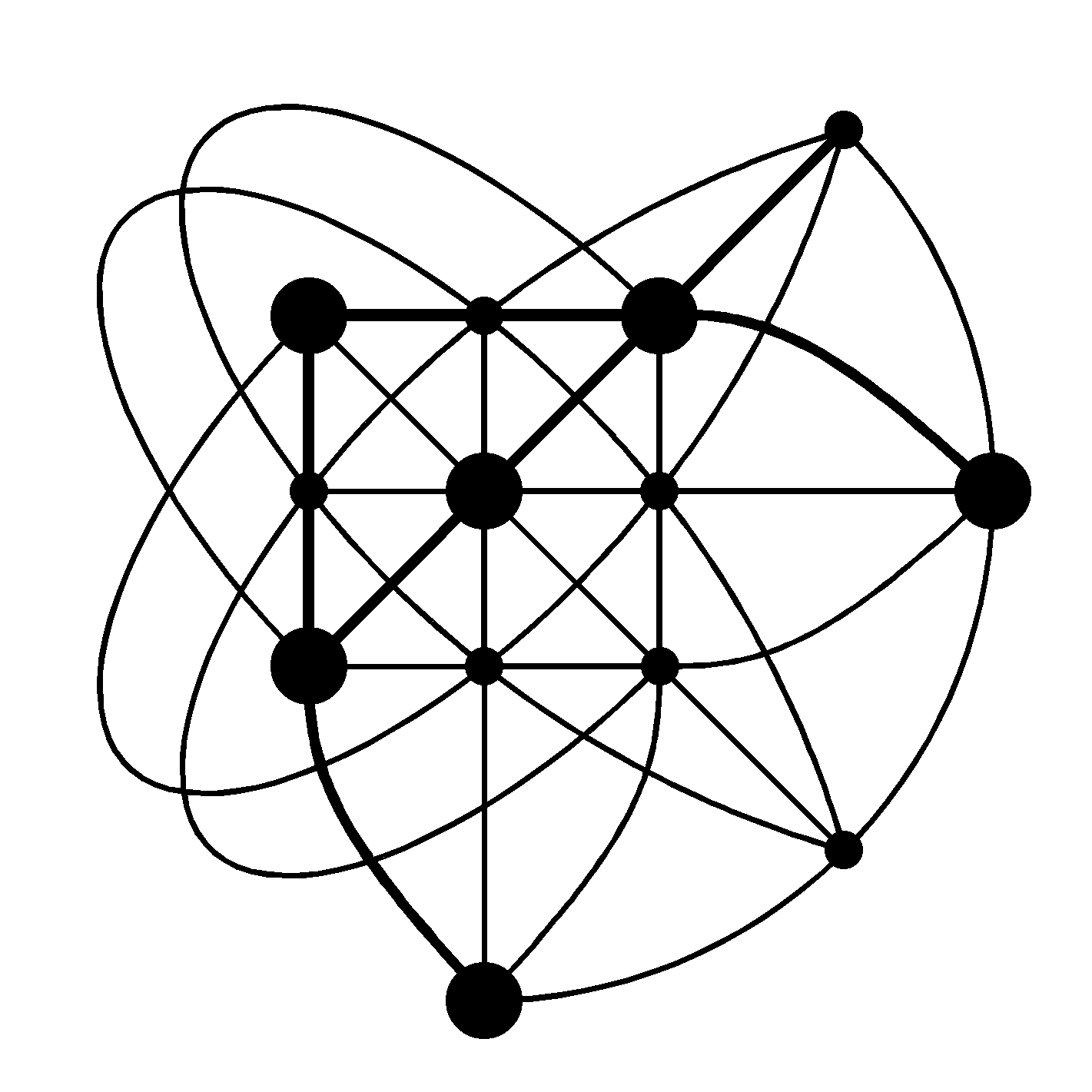}}
        \end{picture}
        {\refstepcounter{abbildung}\label{abb3}
          \centerline{Fig.\ \ref{abb3}.}}
      \end{minipage}
      \begin{minipage}[t]{4.0cm}
         \begin{picture}(4.0,4.0)
         \put(0,0)
         {\includegraphics[width=4.0cm,height=4.0cm]{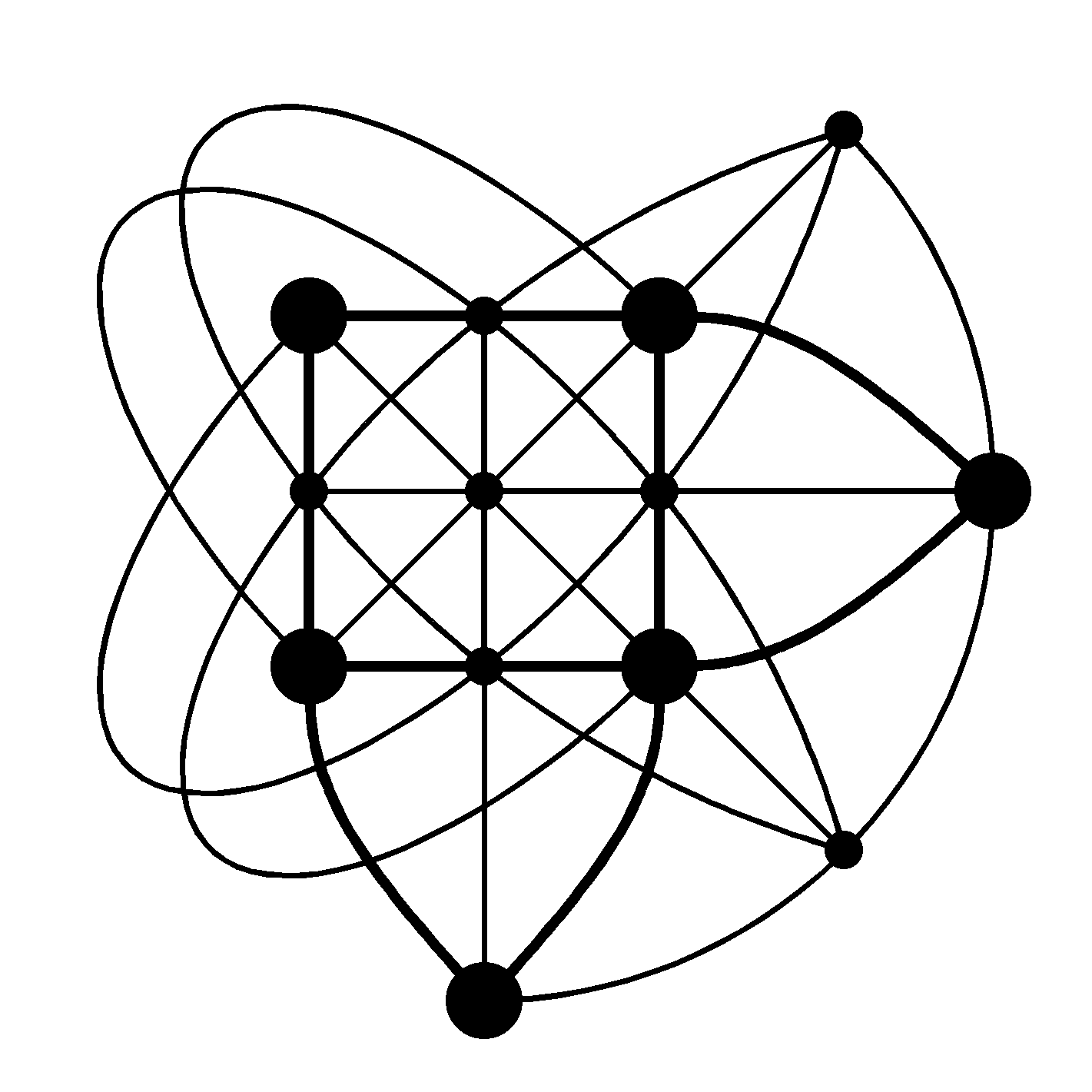}}
        \end{picture}
        {\refstepcounter{abbildung}\label{abb4}
          \centerline{Fig.\ \ref{abb4}.}}
      \end{minipage}
      \end{center}
}%


A set of type $1$ contains a unique line. So there are
exactly
    $13\cdot{9\choose 2} = 13\cdot 36$
sets of type 1.

A set $\cS$ of type $2$ can be written as symmetric
difference of two lines in one way only. Hence there are
exactly
    $\frac{13\cdot 12}{2!} = 13\cdot 6$
sets of type $2$.

If $\cS$ is of type $3$ then each vertex of the ``basic''
triangle is on exactly two trisecants of $\cS$, whereas each
point of the ``inscribed'' triangle is on one trisecant
only. So the role of the two triangles is not the same.
Since two distinct vertices of the inscribed triangle
determine the remaining one uniquely, the number of $6$-sets
of type $3$ is
  $\frac{13\cdot 12 \cdot 9}{3!}\cdot 2\cdot 2 =
  13\cdot 72.$

If $\cS$ is of type $4$ then the defining quadrilateral can
be recovered from $\cS$ as the set its four trisecants. So
the number of sets of type $4$ equals
    $\frac{13\cdot 12\cdot 9\cdot 4}{4!}=13\cdot 18.$

We observe that a $6$-set of type $i\in\{1,2,3,4\}$ has
exactly $i$ trisecants. So the four types of $6$-sets do not
overlap. Finally, from
\begin{displaymath}
  13\cdot(36+6+72+18)=13\cdot 132={13 \choose 6},
\end{displaymath}
our list from above comprises all $6$-sets of points.

Let $U\in\cP$ be a fixed point and put
$\cW:=\cP\setminus\{U\}$. A {\em block}, say $\cB$, is
defined to be a subset of $\cW$ satisfying one of the
following conditions:

\begin{enumerate}
  \item[A.]
$\cB$ is the symmetric difference of two distinct lines,
each not incident with $U$.
  \item[B.]
$\cB\cup\{U\}$ is the union of two distinct lines.
  \item[C.]
$\cB$ consists of a quadrangle together with two of its
diagonal points; moreover, $U$ is the remaining diagonal
point.
\end{enumerate}

If a block $\cB$ is of type A, B, or C, then
$\cP\setminus(\cB\cup\{U\})$ is easily seen to be a $6$-set
of type $1$, $2$ or $4$, respectively. Thus the blocks fall
into classes A, B, and C. Also, let us remark that the
complement in $\cW$ of a block of type A, B, or C is a block
of type B, A, or C, respectively.

The number of blocks of type A is equal to the number of
$2$-sets of lines, both not running through $U$. So it is
    $\frac{9\cdot 8}{2} = 36.$

Blocks of type B are of the form $(a\cup b)\setminus\{U\}$
with lines $a\neq b$ and $U\in a\cup b$. Counting the
possibilities for $a$ and $b$, and taking into account
whether $U$ is on both lines or not, shows that there are
precisely
    $4\cdot 9 + {4\choose 2} = 42$
blocks of type B.

We obtain all quadrangles with diagonal point $U$ by drawing
two distinct lines, say $a$ and $b$, through $U$ and
choosing two distinct points on $a\setminus\{U\}$ and
$b\setminus\{U\}$, respectively. So the number of blocks of
type C equals
    ${4\choose 2}\cdot {3\choose 2}\cdot {3\choose 2} = 54.$

Summing up shows that there are exactly $132$ blocks.

\abstand

Here is our main result:

\begin{thm}
The set $\cW$, together with the set of all blocks, is a
Steiner system $S(5,6,12)$.
\end{thm}

\pf (a) By definition, all blocks have exactly $6$ elements
and $\#\cW=12$.

(b) We show that each $5$-set $\cM$ in $\cW$ belongs to at
least one block. There are four cases, depending on the type
of the $6$-set $\cS:=\cM\cup\{U\}$:
  {\unitlength1.3cm
    \begin{center}
     \begin{minipage}[t]{4.0\unitlength}
      \begin{picture}(4.0,4.0)
         \put(0,0)
         {\includegraphics[width=4.0\unitlength,height=4.0\unitlength]{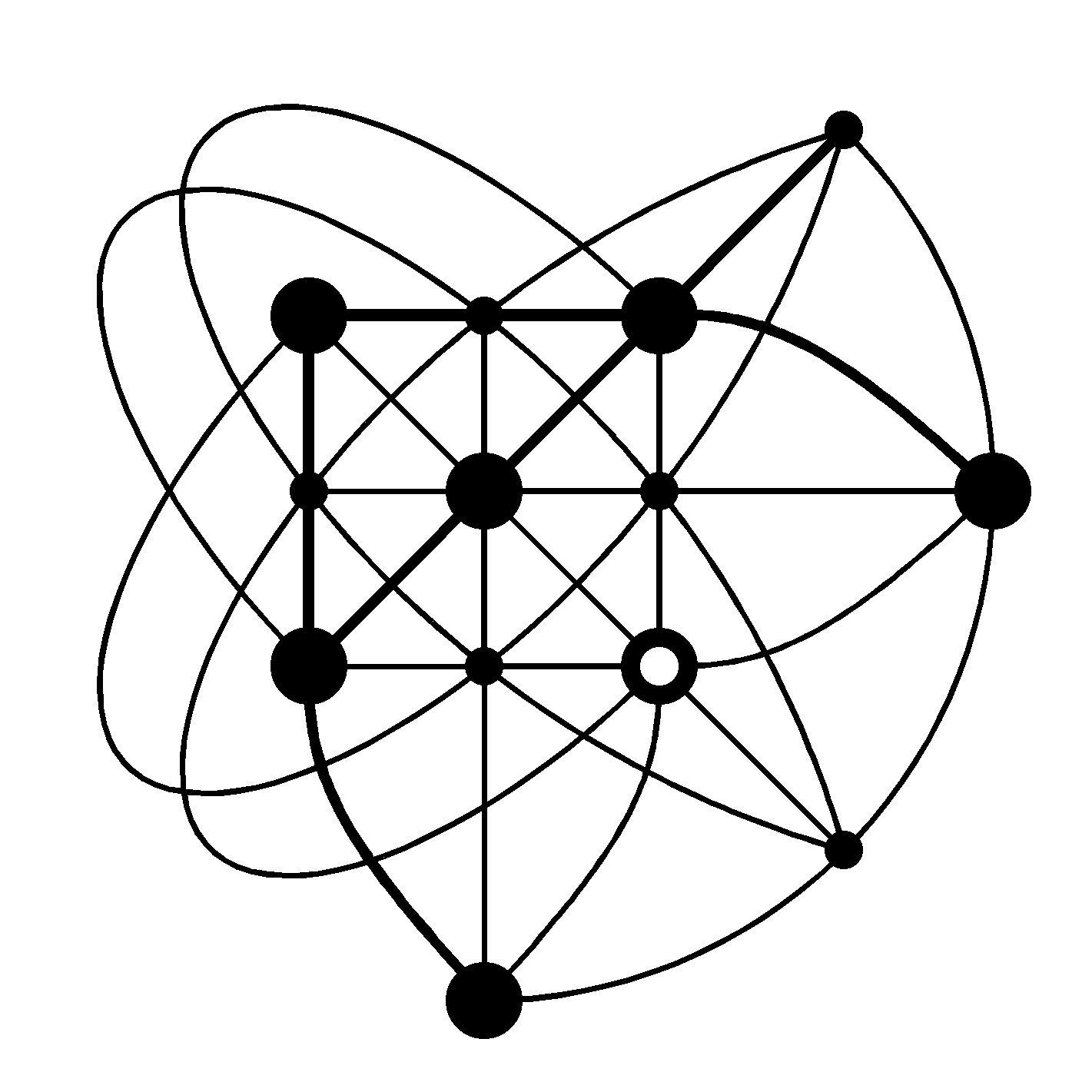}}
         \put(0.7 ,1.22){$A$}
         \put(0.95,3.1){$B$}
         \put(2.25 ,3.1){$C$}
         \put(3.8 ,2.05){$P$}
         \put(1.3 ,2.3){$Q$}
         \put(0.35 ,0.2){$R=U$}
         \put(2.47,1.22){$X$}
      \end{picture}
      {\refstepcounter{abbildung}\label{abb5}
       \centerline{Fig.\ \ref{abb5}.}}
     \end{minipage}
     \begin{minipage}[t]{4.0\unitlength}
      \begin{picture}(4.0,4.0)
         \put(0,0)
         {\includegraphics[width=4.0\unitlength,height=4.0\unitlength]{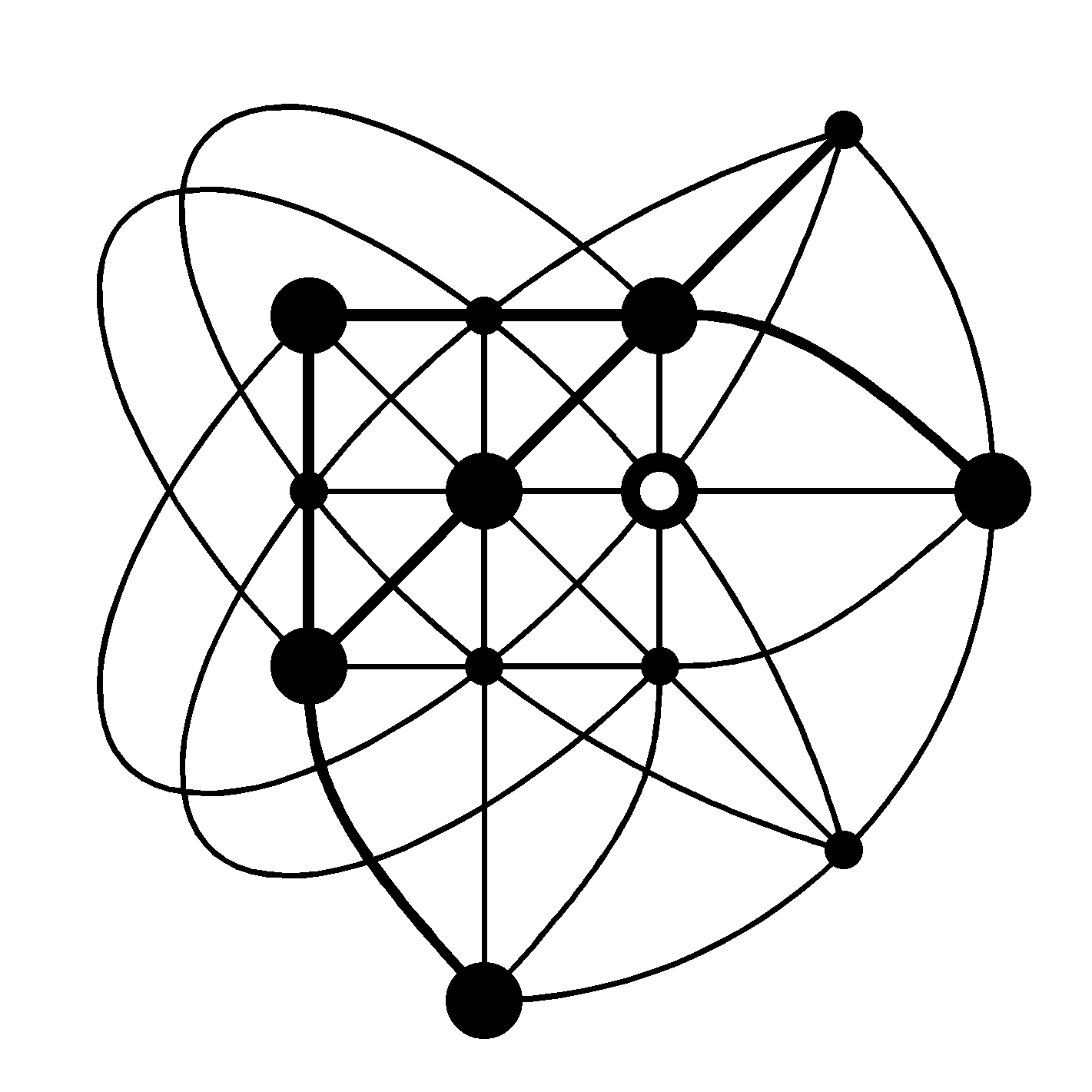}}
         \put(0.7 ,1.22){$A$}
         \put(0.95,3.1){$B$}
         \put(1.5 ,3.1){$C=U$}
         \put(3.8 ,2.05){$P$}
         \put(1.3 ,2.3){$Q$}
         \put(1.25 ,0.2){$R$}
         \put(2.6,2.3){$X$}
      \end{picture}
      {\refstepcounter{abbildung}\label{abb6}
          \centerline{Fig.\ \ref{abb6}.}}
      \end{minipage}
      \begin{minipage}[t]{4.0\unitlength}
      \begin{picture}(4.0,4.0)
         \put(0,0)
         {\includegraphics[width=4.0\unitlength,height=4.0\unitlength]{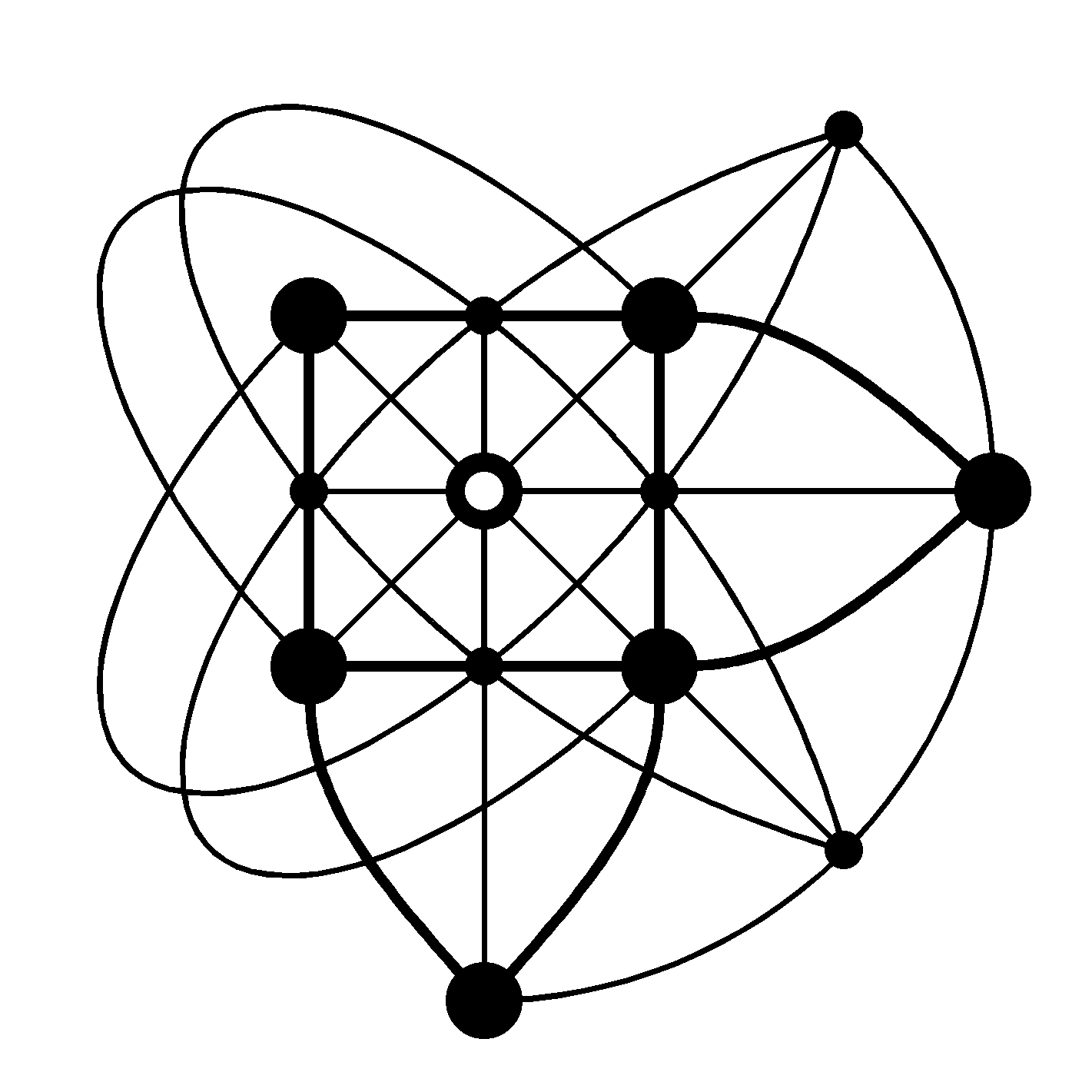}}
         \put(0.7 ,1.22){$A$}
         \put(0.95,3.1){$B$}
         \put(2.25 ,3.1){$C$}
         \put(2.47,1.22){$D$}
         \put(3.8 ,2.05){$E$}
         \put(0.35 ,0.2){$F=U$}
         \put(1.3 ,2.3){$X$}
      \end{picture}
      {\refstepcounter{abbildung}\label{abb7}
          \centerline{Fig.\ \ref{abb7}.}}
     \end{minipage}
    \end{center}
}
\begin{enumerate}
  \item
  Suppose that $\cS$ consists of a line $a$ and two further
  points; let $b$ be the line joining those points. Then $(a\cup b)
  \setminus \{U\}$ is a block of type B containing $\cM$.
  \item
  Let $\cS$ be the symmetric difference of distinct lines
  $a$ and $b$. Then $(a\cup b)\setminus\{U\}$ is a block of
  type B with the required property.
  \item
  Let $\cS$ be the union of a triangle $\{A,B,C\}$ and an
  inscribed triangle $\{P,Q,R\}$ such that $P\in BC$, $Q\in
  CA$, and $R\in AB$.
  There are two subcases:

  If $U\in\{P,Q,R\}$, say $R=U$, then put $\{X\}:=AP\cap
  BQ$. Then $\{X\}\in CR$ and $\{A,B,C,X \}$ is a quadrangle
  with diagonal points $P$, $Q$, and $R=U$ which gives
  rise to a block of type C containing $\cM$ (fig.\ \ref{abb5}).

  If $U\in\{A,B,C\}$, say $C=U$, then put $\{X\}:=PQ\cap RU$.
  Then $U\notin AB \cup PQ$. So the symmetric difference
  of $AB$ and $PQ$ is a block of type A through $\cM$
  (fig.\ \ref{abb6}).
  \item
  Let $\cS=\{A,B,C,D,E,F\}$ be the set of vertices of a
  quadrangle. W.l.o.g.\ let $\{U\}=\{F\}=AB\cap CD$. So
  $\{A,B,C,D\}$ is a quadrangle with diagonal points $E$,
  $F=U$, and $X$, say. Therefore
  $\{A,B,C,D,E,X\}\supset\cM$ is a block of type C
  (fig.\ \ref{abb7}).
\end{enumerate}

(c) Given a $5$-set $\cM\subset\cW$ then denote by $r(\cM)$
the number of blocks passing through it. Since each of the
$132$ blocks contains exactly $6$ subsets of $\cW$ with $5$
elements, we obtain from the principle of counting in two
ways that
\begin{displaymath}
\sum_{\cM\subset\cW, \atop \#\cM=5} { r(\cM)} = 132\cdot 6 =
  792.
\end{displaymath}
From (b), $r(\cM)\geq 1$ for each of the ${12\choose 5}=792$
sets $\cM$ appearing in the sum above. So $r(\cM)= 1$ is
constant. This completes the proof. \qed



{\small Hans Havlicek, Institut f\"ur Geometrie, Technische
Universit\"at, Wiedner Hauptstra{\ss}e 8--10, A--1040 Wien,
Austria.}

\abstand

{\small Hanfried    Lenz,    Mathematisches Institut II
(WE2), Freie Universit\"at Berlin,    Arnimallee 3, D-14195
Berlin, Germany.}

\end{document}